\RequirePackage{fix-cm}
\documentclass[final,times,3p,twocolumn]{elsarticlehack}

\usepackage{amsmath, amsfonts, amsthm, amssymb}
\usepackage{placeins}
\usepackage{pdfpages}
\usepackage[hypertexnames=false,colorlinks=true,breaklinks=true,bookmarks=true,urlcolor=blue,citecolor=blue,linkcolor=blue,bookmarksopen=false,draft=false]{hyperref}
\usepackage[ruled,vlined,linesnumbered,norelsize]{algorithm2e}
\usepackage{comment}

\theoremstyle{plain}
\newtheorem{theorem}{Theorem}

\newtheorem{remark}[theorem]{Remark}

\DeclareMathOperator{\Diag}{Diag}
\DeclareMathOperator{\diag}{diag}
\DeclareMathOperator{\ldet}{ldet}

\makeatletter
\renewcommand*{\top}{%
  {\mathpalette\@transpose{}}%
}
\newcommand*{\@transpose}[2]{%
  \scriptsize
  \raisebox{\depth}{$\m@th#1\mathsf{T}$}%
}
\makeatother

\begin{document}

\begin{frontmatter}



\title{Computing D-Optimal solutions for huge-scale\\ linear and quadratic response-surface models}

%\author[label1,label2,label1]{Gabriel Ponte, Marcia Fampa, Jon Lee}

\author[label1]{Gabriel Ponte}
\author[label2]{Marcia Fampa}
\author[label1]{Jon Lee}

\affiliation[label1]{organization={University of Michigan},
            city={Ann Arbor},
            country={USA}}

\affiliation[label2]{organization={Federal University of Rio de Janeiro},
            country={Brazil}}

% \author{Gabriel Ponte, Marcia Fampa, Jon Lee}

% \affiliation{organization={University of Michigan},%Department and Organization
% %            addressline={}, 
%             city={Ann Arbor},
% %            postcode={}, 
% %            state={},
%             country={USA}}

\begin{abstract}
We consider algorithmic approaches to the D-optimality problem for 
cases where the input design matrix is large and highly structured,
in particular implicitly specified as a full quadratic or linear response-surface model
in several levels of several factors. Using  
row generation techniques of mathematical optimization, in the context of discrete local-search and 
continuous relaxation aimed at branch-and-bound solution, we
are able to design practical  algorithms.

\end{abstract}


% \begin{highlights}
% \end{highlights}


\end{frontmatter}


\section{Introduction}\label{sec:intro}

We consider the D-Optimality problem formulated as 
\begin{equation}\label{prob}\tag{D-Opt}
\begin{array}{lrl}
&\max~ &\ldet\sum_{\ell\in N} x_\ell v_\ell v_\ell^\top\,,\\
&\text{s.t.} & \mathbf{e}^\top x=s,\\
&&0\leq x ,~ x\in\mathbb{Z}^n,
\end{array}
\end{equation}
% \begin{align*}\label{prob}\tag{D-Opt}
% \textstyle
% % &\max \left\{ \ldet\sum_{\ell\in S}  v_\ell v_\ell^\top \, : \, |S|=s,~ S\subset N\right\},\\
% % &\qquad =
% &
% %z:=
% \max \left\{ \ldet\sum_{\ell\in N} x_\ell v_\ell v_\ell^\top \, : \, \mathbf{e}^\top x=s,~  l\leq x \leq u,~ x\in\mathbb{Z}^n
% %x\in \{0,1\}^n
% \right\}.
% % ,\\[4pt]
% % &\qquad=\max \left\{ \ldet \left(
% % \sum_{\ell\in N} l_\ell v_\ell v_\ell^\top 
% % + \sum_{\ell\in N} x_\ell v_\ell v_\ell^\top\right) \, : \right.\\[3pt]
% % &\qquad\qquad\qquad\qquad\, \left.\mathbf{e}^\top x=s-\mathbf{e}^\top l\,,~  0\leq x \leq u-l,~ x\in\mathbb{Z}^n
% % \vphantom{\sum_{\ell\in N} l_\ell v_\ell v_\ell^\top}\right\},
% \end{align*}
where $\ldet$ is the natural logarithm of the determinant,
$v_\ell \in \mathbb{R}^{m}$, 
for $\ell\in N:=\{1,\ldots,n\}$,
%$0 < u\in\mathbb{Z}^n$, 
with
$m\leq s <n$ natural numbers. 

% satisfying 
% $ s \leq \mathbf{e}^\top u$. 

\ref{prob} is a fundamental problem in ``experimental design''.  
Let $A:= (v_1, v_2, \dots, v_n)^\top$, so  
$
\textstyle
\sum_{\ell\in N} x_\ell v_\ell v_\ell^\top = A^\top \Diag(x) A.
$
% We consider the least-squares  
% problem $\min_{\theta\in \mathbb{R}^{m}} \|A\theta -b\|_2$\,, where $b$ is an 
% arbitrary response vector. 
We consider the situation where each $v_\ell$
corresponds to a costly experiment, which could be carried out with repetition, within limits,
and we have a budget to carry out  $s$ experiments.
For a given feasible solution $\tilde{x}$, let $A_{\tilde{x}}$
be a matrix with $v_\ell^\top$ repeated $\tilde x_\ell$ times, for $\ell\in N$, as its rows, and let 
$b_{\tilde{x}}$
denote the associated response vector.
We have the associated least-squares 
problem $\min_{\theta\in \mathbb{R}^m} \|A_{\tilde x}\theta -b_{\tilde{x}}\|_2$\,.
The generalized variance of the least-squares parameter estimator $\hat \theta$
is inversely proportional to $\det \sum_{\ell\in N}  \tilde{x}_\ell v_\ell v_\ell^\top $\,,
and so \ref{prob} models picking the experiments to 
minimize the generalized variance of  $\hat \theta$
(see \cite{Draper} for more details; also see \cite{Wald,Kiefer,KW,Puk,Fedorov} for history and more information). 

There is a lot of  literature on heuristic 
algorithms for \ref{prob}, based on the usual greedy and local-search ideas;  see, for example the references in \cite[Sec. 1]{KoLeeWayne}. 
Upper bounds and branch-and-bound (B{\&}B) approaches aimed at exact solution of \ref{prob}
were investigated in \cite{KoLeeWayne,KoLeeWayne2,Welch,PonteFampaLeeSBPO22,PonteFampaLeeMPB}.

\medskip

\noindent {\bf Our Contribution.} Using row-generation ideas well known in mathematical optimization, we examine local search and B{\&}B for \ref{prob}, in
a couple of significant situations where the matrix $A$ has many rows but is highly structured. In particular, for full quadratic or linear response-surface models
in several levels of several factors (see \S\ref{sec:quad}). 

\begin{remark}
   Some works (e.g., \cite{PonteFampaLeeMPB}) more generally treat $l\leq x \leq u$ constraints in \ref{prob},
   for $0 \leq l  < u\in\mathbb{Z}^n$. In our version (also employed by \cite{Welch} and \cite{nikolov2015randomized}, for examples), we omit these
   constraints to more cleanly expose our ideas.
\end{remark}

\medskip
\noindent {\bf Notation.}
 % We let $\mathbb{S}^n$  (resp., $\mathbb{S}^n_+$~, $\mathbb{S}^n_{++}$)
 % denote the set of symmetric (resp., positive-semidefinite, positive-definite) matrices of order $n$. 
%  We let $\mathbb{S}^n_+$ (resp., $\mathbb{S}^n_{++}$) denote the set of
%  positive semidefinite (resp., definite) symmetric matrices of order $n$. 
 We let $\Diag(x)$ denote the $n\times n$ diagonal matrix with diagonal elements given by the components of $x\in \mathbb{R}^n$, and $\diag(X)$ denote the $n$-dimensional vector with elements given by the diagonal elements of $X\in\mathbb{R}^{n\times n}$.
We denote an all-ones  vector
by $\mathbf{e}$, and an \hbox{$i$-th} standard unit vector by $\mathbf{e}_i$\,.
% and an identity matrix by $I$. 

% For convenience, we let $A:= (v_1, v_2, \dots, v_n)^\top$, so we have  
% $
% \textstyle
% \sum_{\ell\in N} x_\ell v_\ell v_\ell^\top = A^\top \Diag(x) A.
% $
% For matrices $A$ and $B$, 
% $A\bullet B:=\Trace(A^\top B)$ is the matrix dot-product.
% For matrix $A$, we denote row $i$ by $A_{i\cdot}$ and
% column $j$ by $A_{\cdot j}$~. For $X\in\mathbb{S}^n$, we let $d_i(X)$ be the $i$-th largest element of $\diag(X)$, and  $\lambda_i(X)$ be the $i$-th largest eigenvalue of $X$.

\section{Quadratic and linear response-surface models}\label{sec:quad}

A very interesting and classical family of D-optimality problems 
relates to the (full) ``quadratic response-surface model'' of \cite[Section 3]{Welch}. 
Instances are generated 
depending on parameters $L$ (``levels'') and $F$ (``factors'' or ``explanatory variables''). The $m\times n$ design matrix $A$ relates to 
the model $y\approx \theta_0 
+ \sum_{i=1}^F \alpha_i \theta_i 
+ \sum_{i=1}^F \alpha_i^2 \theta_{ii}
+ \sum_{i<j} \alpha_i \alpha_j \theta_{ij}$\,. 
This kind of full quadratic model features in the seminal paper \cite{Box} on response-surface methodology (also see the references in the early survey article \cite{designsurvey}). 

For this quadratic response-surface model,
we create an $m\times n$ design matrix $A$ with $n:=L^F$  and $m:=1+2F+\binom{F}{2}$. 
Each row of $A$, which  has the form
\begin{equation}\label{row_quad}
v^\top:=(1;~ 
\alpha_1,\ldots,\alpha_F\,;~
\alpha_1^2,\ldots,\alpha_F^2\,;~
\alpha_1 \alpha_2\,,\ldots ,\alpha_{F-1} \alpha_F),
\end{equation}
 is identified by the 
levels in $\{0,1,\ldots,L-1\}$ of the factors $\alpha_1,\ldots,\alpha_F$\,.
Ideally, one would like $L=3$
levels (for a full quadratic model), but this would lead to $n=3^F$ experiments in a ``full factorial design'', which
is quite large if we have 
even a modest number of factors.
Hence our interest in the 
D-opt criterion for picking
a best or good set of experiments to fit a response surface.  
The only instance in \cite{Welch} has $L=F=3$, so $n=3^3=27$, and $m=1+6+3=10$. And
experiments were limited to $s=10,\ldots,20$.
% ,
% $l=0$, and $u=s\mathbf{e}$. 

We also consider the simpler and classical linear response-surface model
$y\approx \theta_0 
+ \sum_{i=1}^F \alpha_i \theta_i 
$\,. This leads to instances 
with again 
$n=L^F$, but now $m=1+F$,
where each row of $A$  has the form
\begin{equation}\label{row_lin}
v^\top:=(1;~ 
\alpha_1,\ldots,\alpha_F).
\end{equation}
Normally, with a linear response-surface model, we take $L=2$.

Generally, for these problems,
we assume  that 
$s$ is modest compared to $n$,
so that we can reasonably hope for a sparse optimal solution (local or global) of \ref{prob}. 

\section{Local search}\label{sec:local}

We are interested in these  quadratic and linear response-surface instances when $n:=L^F$ is huge
 --- too big to want to list the rows of $A$ or iterate over them in a local search. We have a feasible solution $\hat{x}$, but we only store the rows of $A$ that are indexed by the support of $\hat{x}$. We let $B:=\sum_\ell \hat{x}_\ell v_\ell v_\ell^\top$\,. The simplest local-search move brings us to $\hat{x}-\mathbf{e}_j+\mathbf{e}_i$\,,
 for some $j$ such that $\hat{x}_j>0$ and $i\not=j$.
 We do explicitly iterate over $j$ in the support of $\hat{x}$ (which we could easily do in parallel).
 But because we do not explicitly generate and store
 the rows of $A$, we will seek to \emph{generate} a row $v$ so that $\ldet(B-v_jv_j^\top + vv^\top)> \ldet(B)$. Now, for the linear
 model,
 a row of $A$ 
  has the form
\eqref{row_lin},
 and for the 
 quadratic model, a row of $A$ 
  has the form
\eqref{row_quad}, 
with $0\leq \alpha_1,\ldots,\alpha_F \leq L-1$ and integer. So we are led to an optimization formulation for generating a row $v$. Finally, assuming that $M:=B-v_j v_j^\top$ is invertible (which is a genuine but mild restriction), we have 
% (see \eqref{detSM})
\[
\ldet(M+vv^\top) = \log(1+v^\top M^{-1}v) +  \ldet(M).
\]

So, for a best local-search move (considering that we are removing $v_j$),
we want to maximize $v^\top M^{-1}v$ over $v$ defined above (depending on whether we are considering the linear or quadratic model). For the linear response-surface model, this particularizes 
to the quadratic-optimization problem
%
       % \begin{align*} \tag{Sub}\label{sub}
       %  &\max~ \alpha^\top M^{-1}\alpha\\   
       % & \alpha_0 = 1,\\
       % & 0\leq \alpha_\ell\leq L-1,\qquad \ell=1,\ldots,F .
       %  \end{align*}
%
\begin{equation}\tag{Sub}\label{sub}
\begin{array}{lrl}
&\max~ &\alpha^\top M^{-1}\alpha,\\
&\text{s.t.} & \alpha_0 = 1,\\
&&0\leq \alpha_\ell\leq L-1,\qquad \ell=1,\ldots,F .
\end{array}
\end{equation}

It is easy to see that the optimal objective value of \ref{sub}
is at least $\det(B)/\det(M)-1$ (considering the $\alpha$ associated with $v_j$), and any objective 
value above that gives a local move with improving objective-function
value. 
% \gabriel{
% % Isn't it $\det(B)/\det(M) - 1$? 
% %
% % If we consider $\alpha = v_j$, then $\ldet(M + \alpha \alpha^\top) = \ldet(B - v_jv_j^\top + v_j v_j^\top) = \ldet(B)$ so
% % \begin{align*}
% %     &\det(M+\alpha \alpha^\top)  = (1+\alpha^\top M^{-1}\alpha)\det(M)\Leftrightarrow\\
% %     &\det(B)  = (1+\alpha^\top M^{-1}\alpha)\det(M)\Leftrightarrow\\
% %     &\det(B)/\det(M)  = (1+\alpha^\top M^{-1}\alpha)\Leftrightarrow\\
% %     &\det(B)/\det(M) - 1  = \alpha^\top M^{-1}\alpha
% % \end{align*}
% Is this a good inequality to add in the model? For example, having a constraint $\alpha^\top M^{-1}\alpha \geq \det(B)/\det(M) - 1$
% }
% \jon{Probably not as a constraint in the model. First of all, it is not a convex inequality, so it would probably put more work on the solver of Sub. Really, what we have is a know incumbent solution (lower bound) that maybe the
% solver can take advantage of if we can pass it or its value to the solver}

Note that this quadratic-optimization problem
has a convex-maximization objective, 
so we relaxed the explicit integrality, and the optimal solution will be at a vertex of the box constraints. So we have only $2^F$ possible solutions to check, rather than the $L^F$ possible integer solutions.
This means that we will only generate points using the extreme levels, $0$ and $L-1$, but for a linear response-surface model, normally we take $L=2$, and so all levels are extreme.
So, in this case, we really have a ``boolean quadric problem''.
For some modest $F$ (e.g., $F=12$),
we can solve \ref{sub} by enumeration. But for more serious $F$ 
(e.g., $F=20$), we could solve \ref{sub} with say \texttt{Gurobi}, and 
probably it would be wise to aggressively use the solver's heuristics
for generating good feasible solutions.

For the quadratic response-surface model, we are instead led to a quartic \emph{integer}-optimization problem over the integer points in an  $F$-dimensional box, which we could approach with  say \texttt{Baron} or \texttt{SCIP}. 

It is rather easy to find an initial set of $m$ rows giving rise to a feasible solution for \ref{prob}.
For the linear response-surface model, we simply take the $m=1+F$ rows 
of the form \eqref{row_lin} that 
correspond to setting $\alpha$ to be 0 and each of the $F$ choices of $\mathbf{e}_i$\,. 
For the quadratic response-surface model, we take the 
$m=1+2F+\binom{F}{2} $ rows 
of the form \eqref{row_quad} that 
correspond  to taking $\alpha$ to be: (i)  0, (ii)  the $F$ choices of  $\mathbf{e}_i$\,,
(iii)  the $F$ choices of $2\mathbf{e}_i$\,, and (iv)
the $\binom{F}{2}$ choices of $\mathbf{e}_i+ \mathbf{e}_j$\,, for all distinct pairs $i,j$.
In both cases, it is easy to check that the rows form an $m\times m$ invertible matrix, and
so the corresponding solution, taking each of these once, is in
the domain of the objective function of \ref{prob}. 
Moreover, using multiplicities, we can easily see that 
we can get such a solution that is feasible for \ref{prob}. 

\section{Relaxation}\label{sec:relax}

There are a variety of 
upper bounds for the optimal value of \ref{prob} (summarized and developed in \cite{PonteFampaLeeMPB}; and also \cite{KoLeeWayne,KoLeeWayne2,PonteFampaLeeSBPO22}). We focus our attention on one of the best ones, obtained 
as the natural convex continuous relaxation of \ref{prob}, and formulated as 
\begin{equation}\label{cont_rel}\tag{$\mathcal{N}$}
\begin{array}{lrl}
&\max~ &\ldet A^\top \Diag(x) A,\\
&\text{s.t.} & \mathbf{e}^\top x=s,\\
&&0 \leq x,~ x\in\mathbb{R}^n,
\end{array}
\end{equation}
% \begin{equation}\label{cont_rel}\tag{$\mathcal{N}$}
% %z_{\mbox{\tiny $\mathcal{N}$}}
% % \hypertarget{znaturaltarget}{\znaturalthing}:=
% \max \left\{ \ldet \left(A^\top \Diag(x) A\right) \, : \, \mathbf{e}^\top x=s, \, l\leq x\leq u,~ x\in\mathbb{R}^n\right\},
% \end{equation}
with the optimal value referred  to  
% $z_{\mbox{\tiny $\mathcal{N}$}}$ 
as the \emph{natural bound} for \ref{prob}.

We will find it useful to consider the Lagrangian dual  of \ref{cont_rel}, which can be formulated as
% \begin{equation}\label{eq:lag_with_theta}\tag{Du-$\mathcal{N}_{\tiny{\mbox{D-Opt}}}$}
% \begin{array}{lrl}
% &\min &-\ldet \Theta    - \omega^\top l + \nu^\top u + \tau s - {m},\\
% &\text{s.t.} 
% &\Theta \bullet v_iv_i^\top  + \omega_i - \nu_i - \tau = 0,\quad i \in N,\\
% &&\Theta \succ 0,\nu \geq 0, \omega \geq 0.
% \end{array}
% \end{equation}
% \begin{equation}\label{eq:lag_with_theta}\tag{Du-$\mathcal{N}$}
% \begin{array}{lrl}
% &\min &-\ldet \Theta     + \nu^\top u + \tau s - {m},\\
% &\text{s.t.} 
% &\diag(A\Theta A^\top)  + \omega - \nu - \tau\mathbf{e} = 0,\\
% &&\Theta \succ 0,\nu \geq 0, \omega \geq 0.
% \end{array}
% \end{equation}
\begin{equation}\label{eq:lag_with_theta}\tag{Du-$\mathcal{N}$}
\begin{array}{lrl}
&\min &-\ldet \Theta    + \tau s - {m},\\
&\text{s.t.} 
&\Theta \bullet v_iv_i^\top  - \tau \leq 0,\quad i \in N,\\
&&\Theta \succ 0.
\end{array}
\end{equation}
(see \cite{PonteFampaLeeMPB}).

%\subsection{Solving the continuous relaxation of \ref{prob}}

With an eye toward developing a branch-and-bound algorithm
for large-scale instances, aimed at large $n$ and modest $s$, we wish to solve the continuous relaxation of \ref{prob}, without explicitly working with all of the rows of the entire $A$ matrix. 
We assume that we have solved
\ref{cont_rel} and \ref{eq:lag_with_theta} over the restriction defined by a subset of the rows of $A$.
% With our assumption that $l=0$, the dual 
% \ref{eq:lag_with_theta}
% takes the form
% \begin{equation}\label{eq:lag_with_theta2}
% \begin{array}{lrl}
% &\min &-\ldet \Theta    + \nu^\top u + \tau s - {m},\\
% &\text{s.t.} 
% &\Theta \bullet v_iv_i^\top  - \nu_i - \tau \leq 0,\quad i \in N,\\
% &&\Theta \succ 0,\nu \geq 0.
% \end{array}
% \end{equation}
Let $(\hat{\Theta},\hat{\tau})$ be an optimal  solution of the dual \ref{eq:lag_with_theta} (probably generated by a solver applied to the primal \ref{cont_rel}).
We can now formulate the separation problem of maximizing $ \hat{\Theta} \bullet v v^\top$\,,
subject to $v^\top$ being a row of the full $A$.
For the linear response-surface model, this is 
just like \ref{sub}, but with $M^{-1}$ replaced by 
$\hat{\Theta}$. For the quadratic response-surface model, we are instead led to a quartic \emph{integer} optimization problem over a $F$-dimensional box.
In either case, 
we obtain a violated inequality
when the optimal value is greater than $\hat{\tau}$.
From a new $v_i$\,, we get a new row for the restriction.
When we solve the continuous relaxation of the new restriction, some $x_i$ may become 0, and we can consider dropping such $x_i$ from the current restriction. In short, we can apply all of the usual ideas of ``column generation''. We note that we can initialize the row-generation procedure
using exactly the initial solutions that we describe for local search.
In this way, \ref{eq:lag_with_theta} will have an optimal solution, which we can use to 
start the row generation.

% \jon{The difficulty in all of this is that it seems that we really need to let the row generation
% fully converge to get a certifiable upper bound for the optimal value of \eqref{prob}.
% But can we stop before full convergence
% and somehow use the $G(\hat\Theta)$ construction from Section 2 of \cite{PonteFampaLeeMPB}
% without writing the full $A$?}

% %%%%

% \mf{I think we only  get an upper bound when the separation problem does not obtain an inequality violated  by the current dual solution.}

Further, we note that we can  construct a closed-form  feasible solution $(\tilde \Theta, \tilde \tau)$ of \ref{eq:lag_with_theta}, over a restriction defined by a subset $R\subset N$, of the rows  of $A$,  from any feasible solution   $\tilde x$ of \ref{cont_rel} over the same restriction. We assume that   $A_r^\top\Diag(\tilde x)A_r\!\in \!\mathbb{S}^m_{++}$\,, where $A_r$ is the submatrix of $A$ defined by the subset of rows $R$, and we set $\tilde\Theta:= (A_r^\top \Diag(\tilde x) A_r)^{-1}$ and $\hat \tau := \max\{\tilde\Theta \bullet v_iv_i^\top\, :\, i\in R\}$.

\section*{Acknowledgments} 
M. Fampa was supported in part by CNPq grants 305444/2019-0 and 434683/2018-3.  
J. Lee was supported in part by AFOSR grant FA9550-22-1-0172. 
This work is partially based upon work supported by the 
National Science Foundation under Grant No. DMS-1929284 while 
the authors were in residence at the Institute for Computational and Experimental Research in 
Mathematics (ICERM) in Providence, RI, during the Discrete Optimization program.



\bibliographystyle{alpha} 
\bibliography{hugebib}

\begin{thebibliography}{KLW98}

\bibitem[BW51]{Box}
George~E.P. Box and Kenneth~B. Wilson.
\newblock On the experimental attainment of optimum conditions.
\newblock {\em Journal of the Royal Statistical Society: Series B
  (Methodological)}, 13(1):1--38, 1951.

\bibitem[Fed72]{Fedorov}
Valerii~V. Fedorov.
\newblock {\em Theory of Optimal Experiments}.
\newblock Academic Press, New York-London, 1972.
\newblock Translated from the Russian and edited by W. J. Studden and E. M.
  Klimko.

\bibitem[HH66]{designsurvey}
William~J. Hill and William~G. Hunter.
\newblock A review of response surface methodology: A literature survey.
\newblock {\em Technometrics}, 8(4):571--590, 1966.

\bibitem[Kie58]{Kiefer}
Jack Kiefer.
\newblock On the nonrandomized optimality and randomized nonoptimality of
  symmetrical designs.
\newblock {\em The Annals of Mathematical Statistics}, 29(3):675--699, 1958.

\bibitem[KLW94]{KoLeeWayne2}
Chun-Wa Ko, Jon Lee, and Kevin Wayne.
\newblock A spectral bound for {D}-optimality, 1994.
\newblock Unpublished.

\bibitem[KLW98]{KoLeeWayne}
Chun-Wa Ko, Jon Lee, and Kevin Wayne.
\newblock Comparison of spectral and {H}adamard bounds for {D}-optimality.
\newblock In {\em M{ODA} 5}, pages 21--29. Physica, Heidelberg, 1998.

\bibitem[KW59]{KW}
Jack Kiefer and Jacob Wolfowitz.
\newblock Optimum designs in regression problems.
\newblock {\em The Annals of Mathematical Statistics}, 30(2):271--294, 1959.

\bibitem[Nik15]{nikolov2015randomized}
Aleksandar Nikolov.
\newblock Randomized rounding for the largest simplex problem.
\newblock In {\em Proceedings of STOC 2015}, pages 861--870, 2015.

\bibitem[PFL22]{PonteFampaLeeSBPO22}
Gabriel Ponte, Marcia Fampa, and Jon Lee.
\newblock Exact and heuristic solution approaches for the {D}-optimality
  problem, 2022.
\newblock In: Proceedings of the LIV Brazilian Symposium on Operations Research
  (SOBRAPO 2022),
  \url{https://proceedings.science/proceedings/100311/_papers/157447/download/fulltext_file2}.

\bibitem[PFL23]{PonteFampaLeeMPB}
Gabriel Ponte, Marcia Fampa, and Jon Lee.
\newblock Branch-and-bound for {D}-optimality with fast local search and
  variable-bound tightening, 2023.
\newblock preprint.

\bibitem[Puk06]{Puk}
Friedrich Pukelsheim.
\newblock {\em Optimal Design of Experiments}, volume~50 of {\em Classics in
  Applied Mathematics}.
\newblock SIAM, 2006.
\newblock Reprint of the 1993 original.

\bibitem[SJD75]{Draper}
Ralph~C. St.~John and Norman~R. Draper.
\newblock D-optimality for regression designs: A review.
\newblock {\em Technometrics}, 17(1):15--23, 1975.

\bibitem[Wal43]{Wald}
Abraham Wald.
\newblock On the efficient design of statistical investigations.
\newblock {\em Annals of Mathematical Statistics}, 14:134--140, 1943.

\bibitem[Wel82]{Welch}
William~J. Welch.
\newblock Branch-and-bound search for experimental designs based on
  {D}-optimality and other criteria.
\newblock {\em Technometrics}, 24(1):41--48, 1982.

\end{thebibliography}





% \end{thebibliography}
\end{document}